\documentclass[11pt]{article}
\usepackage{latexsym}
\usepackage[all]{xy}
\usepackage{amsmath}
\usepackage{amssymb}
\usepackage{amsfonts}

\newtheorem{thm}{Theorem}[subsection]
\newtheorem{prop}[thm]{Proposition}

\newtheorem{df}[thm]{Definition}

\newtheorem{rmk}[thm]{Remark}
\newtheorem{ex}[thm]{Example}
\newtheorem{Q}[thm]{Question}

\begin{document}

\title{\textbf{A sketchy note on \\enriched homotopical topologies \\and enriched homotopical stacks }} 
\bigskip
\bigskip

\author{
Gabriele Vezzosi \\
\small{Dipartimento di Matematica Applicata}\\
\small{Universit\`a di Firenze}\\
\small{Italy}\\
\bigskip}

\date{\textit{}}
\maketitle

\tableofcontents


\begin{section}{Introduction}
The purpose of this rough note (whose bulk was written two years ago, and is not intended for publication) is merely to begin the exploration of the theory of stacks over enriched homotopical sites. The word exploration here is essential: one sets up the basics of the theory with some examples in mind while the rest is uncertain.\\ An enriched homotopical site is a category $C$ enriched over a symmetric monoidal model category $M$ together with a notion of topology (an enriched homotopy topology, for lack of a better name) which is defined through the enrichment and is compatible, in some sense, with the model structure on $M$. The example to keep in mind (though the corresponding theory of stacks already exists, see \cite{hagI}) is when $M=\mathbf{SSet}$, i.e. when $C$ is an $S$-category. In this case stacks are particular simplicial presheaves, i.e. $\mathbf{SSet}$-functors from $C^{op}$ to $\mathbf{SSet}$. In the general case, stacks will be particular $M$-functors from $C^{op}$ to $M$. Unlike in \cite{hagI} (for the case $M=\mathbf{SSet}$) I have tried to keep the definition of an enriched homotopy topology in the sieves' style instead of the covering families' one (see however Section 5), but this might be considered as a minor issue.\\
Unfortunately the theory sketched in this note is seriously limited by the fact  that it does not seem to be the right one if the enrichment takes place in a \textit{stable} model category $M$; this is essentially due to the fact that the notion of homotopy monomorphism (used below to define enriched homotopy sieves) is in this case equivalent to the notion of weak equivalence. One would really like to have e.g. a new and useful notion of enriched homotopy topology for objects like $\mathrm{B}A$, where $A$ is a ring spectrum; this note sheds no interesting light on this problem\footnote{In Section 4.4 and 5, however, we have tried to hint at an alternative route that might solve this problem.}. Therefore the theory sketched below will be possibly useful only in the case of enrichment in non-stable model categories, e.g. for higher categories.\\
As shown in \cite[Prop. 4.5]{RSS}, any (proper, cofibrantly generated) stable model category is Quillen equivalent to a (canonically defined) simplicial model category. It is natural to ask whether it is possible to extend this result to the monoidal case, and further to enrichments.
On the other hand, I prefer to think of (weak) enrichments in symmetric monoidal Segal or simplicial categories (\cite{totann}) as the correct theory of ``homotopical'' enrichments, i.e. of enrichments in a category that carries suitably compatible monoidal and homotopical structures (of which symmetric monoidal model categories are only one possible model). From this point of view, a theory of weakly enriched homotopy topologies (and stacks) would be much more interesting than the tentative theory sketched below, but I guess one needs new ideas to deal with this question properly. \\

As a note, the treatment is everywhere sketchy and somewhere dry; some proofs are deferred, since they do not seem to add any clearness to the conceptual picture which was my main concern here. Many theoretical points remained unexplored (e.g. enriched homotopy topoi) and most of all, many possibly interesting applications have been only suggested as questions to the interested reader. I hope somebody will deal with these issues in the future, just to see how fruitful (or not) they might be.\\

\textbf{Acknowledgments.} This note was very much shaped by discussions with Bertrand To\"en as a follow-up to \cite{hagI}; I thank him sincerely. I am indebted to Charles Rezk for his notion of homotopy monomorphism, to Andr\'e Joyal and Joseph Tapia for discussions and suggestions. \\
\end{section}

\begin{section}{Enriched homotopical topologies}
Let $M$ be a symmetric monoidal model category which is cofibrantly generated (as a model category). Let $C$ be a category enriched over $M$ that for technical reason we will suppose having cofibrant $M$-Hom's; we will write $C_{0}$ for its underlying (usual) category. Viewing $M$ with its canonical $M$-enrichment, let us denote by $\mathsf{\mathbf{Pr}}_{M}(C)$ the category of $M$-functors $C^{op}\rightarrow M$; this is a model category enriched over $M$ with the so-called projective model structure (\cite[Thm. 13.10.17]{hi}). \\

Since $C$ is $M$-enriched, for any object $x$ in $C$, we have an object $\underline{h}_{x} \in \mathsf{\mathbf{Pr}}_{M}(C)$ defined as $y \mapsto \underline{Hom}_{C}(y,x)$, where $\underline{Hom}_{C}(y,x)$ denotes the Hom-object in $M$ of morphisms from $x$ to $y$ in $C$. 

\begin{df}\label{hocat}
If $C$ and $M$ are as above we define the \emph{homotopy category} of $C$ \emph{with respect to} $M$ as the category $\mathsf{\mathbf{Ho}}_{M}(C)$ with the same objects as $C$ and with $$Hom_{\mathsf{\mathbf{Ho}}_{M}(C)}(x,y):=Hom_{\mathsf{\mathbf{Ho}}(\mathsf{\mathbf{Pr}}_{M}(C))}(\underline{h}_{x},\underline{h}_{y})$$ for any pair $(x,y)$ of objects in $C$.
\end{df}

\begin{rmk}
\emph{When $M=\mathsf{\mathbf{SSet}}$, i.e. $C$ is an $S$-category, then $\mathsf{\mathbf{Ho}}_{\mathsf{\mathbf{SSet}}}(C)$ is equivalent to the category $\mathsf{\mathbf{\pi}}_{0}C$ of connected components of $C$. This follows from the fact that by the $\mathsf{\mathbf{SSet}}$-enriched Yoneda lemma, $\underline{h}_{x}$ is cofibrant in $\mathsf{\mathbf{Pr}}_{\mathsf{\mathbf{SSet}}}(C)$ and from the fact that there is an isomorphism $$Hom_{\mathsf{\mathbf{Ho}}(\mathsf{\mathbf{Pr}}_{\mathsf{\mathbf{SSet}}}(C))}(F,G)\simeq \pi_{0}(\mathbb{R}\underline{Hom}_{\mathsf{\mathbf{Pr}}_{\mathsf{\mathbf{SSet}}}(C)}(F,G)),$$ for any $F$ and $G$ in $\mathsf{\mathbf{Pr}}_{\mathsf{\mathbf{SSet}}}(C)$, since $\mathsf{\mathbf{Pr}}_{\mathsf{\mathbf{SSet}}}(C)$ is a simplicial model category.}
\end{rmk}

\noindent Recall the following definition of homotopy monomorphisms due to C. Rezk.

\begin{df} Let $N$ be a right proper model category. A morphism $f:x\rightarrow y$ in $\mathsf{\mathbf{Ho}}(N)$ is called a \emph{homotopy monomorphism} if given any morphism $f':x'\rightarrow y'$ in $N$ such that $\mathsf\mathbf{{Ho}}(f')$ is isomorphic to $f$ in the category of arrows in $\mathsf{\mathrm{Ho}}(N)$, the following square is homotopy cartesian (\cite[13.3.11]{hi})

$$\xymatrix{ x \ar[r]^{\mathrm{id}} \ar[d]_{\mathrm{id}} & x\ar[d]^-{f'} \\ 
             x \ar[r]_{f'}                                & y }$$

\end{df}
This definition makes sense (i.e. it is independent of the choice of the lifting $f'$ of $f$ to $N$) due to the invariance of homotopy cartesian squares with respect to objectwise equivalences (\cite[Prop. 13.3.13]{hi}).\\
Note that if the model structure is trivial then a map is an homotopy monomorphism iff it is a monomorphism in the usual categorical sense. There is a by now obvious dual notion of \textit{homotopy epimorphism}.\\

\begin{prop}\label{hpbmono}
The homotopy pullback of a homotopy monomorphism is a homotopy monomorphism.
\end{prop}

\textbf{Proof.} Let $f: x\rightarrow y$ be our homotopy mono and $g: y'\rightarrow y$ be a map in $N$; we may suppose that $f$ is already lifted to $N$ so let us factor it as $x\rightarrow x'\rightarrow y$, a trivial cofibration composed with a fibration. The homotopy pull-back of $f$ along $g$ is then represented by the pull-back $\widetilde{y'}:=x'\times_{y}y'\rightarrow y'$ (\cite[Cor. 13.3.8]{hi}) and we must prove that the canonical map $$\widetilde{y'}\longrightarrow\widetilde{y'}\times_{y'}\widetilde{y'}$$ is an equivalence. Since the diagonal map $\Delta: x'\rightarrow x'\times_{y}x'$ is an equivalence (by hypohesis) and $ x'\times_{y}x'\rightarrow y$ is a fibration, $\widetilde{y'}\longrightarrow\widetilde{y'}\times_{y'}\widetilde{y'}$ is an equivalence by \cite[Prop. 13.3.9]{hi}. \hfill $\mathbf{\Box}$

\begin{prop}\label{monoprop}
\begin{itemize}
	\item Any equivalence in $N$ is a homotopy monomorphism.
	\item The class of homotopy monomorphisms is closed under composition.
\end{itemize}
\end{prop}
\textbf{Proof.} Easy verification. \hfill $\mathbf{\Box}$

\begin{ex}\label{cat}
\emph{A nice and useful example is the case $N=\mathsf{\mathbf{Cat}}$ with the model structure described in \cite{re}, where equivalences are categorical equivalences and cofibrations are functors which are injective on objects. In this case, a functor is an homotopy monomorphism iff it is faithful and surjective (hence bijective) on isomorphisms. So, fully faithful functors are homotopy monomorphisms but the class of homotopy monomorphisms is obviously strictly bigger; for example, if $k$ is a ring, the obvious functor from the category of complexes of $k$-modules with morphisms given by quasi-isomorphisms to the category of complexes of $k$-modules (with all morphisms) is an homotopy monomorphism but is not fully faithful.}
\end{ex}

\begin{rmk}\label{mono}
\emph{Let $M=\mathbf{SSet}$ and $C$ be a usual category, endowed with the trivial $\mathbf{SSet}$-enrichment. Then $\mathbf{Pr}_{\mathbf{SSet}}(C)=\mathbf{SPr}(C)$, the category of simplicial presheaves on $C$. If $F$ denotes a presheaf of sets on $C$, we may view it as a constant simplicial presheaf on $C$; then, a morphism $F\rightarrow G$ in $\mathbf{SPr}(C)$ between constant simplicial presheaves is a homotopy monomorphism iff it is a monomorphism.} 
\end{rmk}

\begin{df}\label{Msieve} If $x$ is an object in $C$, a \textit{homotopy} $M$\emph{-sieve} $\underline{R}$ over $x$ in $C$ is a homotopy monomorphism $\underline{R} \rightarrow \underline{h}_{x}$ in the homotopy category $\mathsf{\mathbf{Ho}}(\mathsf{\mathbf{Pr}}_{M}(C))$. 
\end{df}

\begin{rmk}\label{sievesasehsieves}
\emph{Let $M=\mathbf{SSet}$ and $C$ be a usual category, endowed with the trivial $\mathbf{SSet}$-enrichment; then, as observed above, $\mathbf{Pr}_{\mathbf{SSet}}(C)=\mathbf{SPr}(C)$. Since the enrichment is trivial (i.e., constant), by Remark \ref{mono}, any standard sieve in $C$ is a  $\mathbf{SSet}$-sieve in $C$. Viceversa, any $\mathbf{SSet}$-sieve $F\rightarrow \underline{h}_{x}$ over $x$ in $C$, induces (by adjunction) a usual sieve $\pi_{0}F \rightarrow h_{x}$ in $C$. } 
\end{rmk}

\begin{ex} \emph{Consider a (symmetric) monoidal model category $\mathcal{M}$ and view it as a $\mathsf{\mathbf{Cat}}$-enriched category with one single object $*$ whose endomorphism category is $\mathcal{M}$ itself. Then the category $\mathsf{\mathbf{Pr}}_{\mathsf{\mathbf{Cat}}}(\mathcal{M})$ is naturally equivalent to the category $\mathsf{\mathbf{Mod}}_{\mathcal{M}}$ of $\mathcal{M}$-Modules (i.e categories with an action of $\mathcal{M}$); moreover, under this equivalence, the model structure on $\mathsf{\mathbf{Pr}}_{\mathsf{\mathbf{Cat}}}(\mathcal{M})$ translates into the model structure on $\mathsf{\mathbf{Mod}}_{\mathcal{M}}$ created by the forgetful functor $\varphi: \mathsf{\mathbf{Mod}}_{\mathcal{M}}\longrightarrow \mathsf{\mathbf{Cat}}$ (i.e. a morphism in $\mathsf{\mathbf{Mod}}_{\mathcal{M}}$ is a fibration/equivalence/cofibration iff the same is true for its image under $\varphi$, for the model structure in $\mathsf{\mathbf{Cat}}$ described in Example \ref{cat}). Therefore, a homotopy  $\mathsf{\mathbf{Cat}}$-sieve in $\mathsf{\mathbf{Pr}}_{\mathsf{\mathbf{Cat}}}(\mathcal{M})$ is identified with a morphism $R\rightarrow \mathcal{M}$ of $\mathcal{M}$-Modules which, as a functor, is faithful and surjective on isomorphisms.}
\end{ex}

If $\underline{R}\rightarrow \underline{h}_{x}$ is a homotopy $M$-sieve for $C$ and $f:y\rightarrow x$ is a morphism in $C_{0}$, there is an \textit{inverse image} homotopy $M$-sieve $f^{-1}\underline{R}\rightarrow \underline{h}_{y}$ defined by the following homotopy pullback square 

$$\xymatrix{ f^{-1}\underline{R} \ar[r] \ar[d] & \underline{R} \ar[d] \\ 
             \underline{h}_{y} \ar[r]          & \underline{h}_{x} }$$

Here we used the weak version of the $M$-enriched Yoneda lemma: $Hom_{\mathsf{\mathbf{Pr}}_{M}(C)}(\underline{h}_{z},F)=Hom_{M}(1,F(z))$ (bijection of sets).

\begin{df}
A morphism $f:y\rightarrow x$ in $C$ is said \emph{to belong} to an homotopy $M$-sieve $\underline{R}\rightarrow \underline{h}_{x}$ over $x$ if the induced morphism $\mathsf{\mathbf{Ho}}(f:\underline{h}_{y}\rightarrow\underline{h}_{x})$ factors through $\underline{R}\rightarrow \underline{h}_{x}$ in $\mathsf{\mathbf{Ho}}(\mathsf{\mathbf{Pr}}_{M}(C))$, i.e. there is a commutative diagram in $\mathsf{\mathbf{Ho}}(\mathsf{\mathbf{Pr}}_{M}(C))$ $$\xymatrix{\underline{R} \ar[rr] & & \underline{h}_{x}\\
& & \underline{h}_{y}. \ar[u]_-{\mathsf\mathbf{{Ho}}(f)} \ar[ull]
}$$
\end{df}

\begin{df}\label{Mtop}
An $M$\emph{-enriched homotopy topology} $\underline{\tau}$ on $C$ consists of data $\underline{Cov}_{\underline{\tau}}(x)$ of sets of homotopy $M$-sieves on $x$ in $C$, for any object $x$ in $C$, satisfying the following conditions:
\begin{enumerate}
  \item for any $x$ in $C$, $\underline{Cov}_{\underline{\tau}}(x)$ is stable under equivalences, i.e. if $u:\underline{R}\rightarrow \underline{h}_{x}$ belongs to $\underline{Cov}_{\underline{\tau}}(x)$ and $u':\underline{R'}\rightarrow \underline{h}_{x}$ is an $M$-enriched homotopy sieve over $x$ such that $u\simeq u'$ in $\mathsf\mathbf{{Ho}}(\mathsf{\mathbf{Pr}}_{M}(C))/\underline{h}_{x}$, then $u'$ belongs to $\underline{Cov}_{\underline{\tau}}(x)$;
	\item  for any $x$ in $C$, the homotopy sieve $(\mathrm{id}:\underline{h}_{x}\rightarrow \underline{h}_{x})$ belongs to $\underline{Cov}_{\underline{\tau}}(x)$;
	\item if $\underline{R}\rightarrow \underline{h}_{x}$ belongs to $\underline{Cov}_{\underline{\tau}}(x)$ and $f:y\rightarrow x$ is a morphism in $C_{0}$, then $f^{-1}\underline{R}\rightarrow \underline{h}_{y}$ belongs to $\underline{Cov}_{\underline{\tau}}(y)$;
	\item if $\underline{R}\rightarrow \underline{h}_{x}$ belongs to $\underline{Cov}_{\underline{\tau}}(x)$ and $\underline{S}\rightarrow \underline{h}_{x}$ is an homotopy $M$-sieve over $x$ such that $f^{-1}\underline{R}\rightarrow \underline{h}_{y}$ belongs to $\underline{Cov}_{\underline{\tau}}(y)$, for any $f:y\rightarrow x$ belonging to $\underline{R}\rightarrow \underline{h}_{x}$, then $\underline{S}\rightarrow \underline{h}_{x}$ belongs to $\underline{Cov}_{\underline{\tau}}(x)$.
\end{enumerate}
\noindent The pair $(C,\underline{\tau})$ will be called an $M$\emph{-enriched homotopy site}.
\end{df}

\begin{rmk}\label{trivial}
\emph{Let $M=\mathbf{SSet}$ and $C$ be a usual category, endowed with the trivial $\mathbf{SSet}$-enrichment; then, as observed above, $\mathrm{Pr}_{\mathbf{SSet}}(C)=\mathbf{SPr}(C)$. Since the enrichment is trivial (i.e., constant), by Remark \ref{sievesasehsieves}, any Grothendieck topology on $C$ gives  a  $\mathbf{SSet}$-enriched homotopy topology on $C$.} 
\end{rmk}

Let $\underline{\tau}$ be an enriched homotopy topology  on an $M$-enriched category $C$ and $\underline{R}\rightarrow \underline{h}_{x}$ be an enriched $\underline{\tau}$-covering sieve. Define  $\mathsf{\mathbf{Ho}}(\underline{R})$ as the set of maps $f:y\rightarrow x$ in $\mathsf{\mathbf{Ho}}_{M}(C)$ such that there is a factorization $$\xymatrix{\underline{R} \ar[r] & \underline{h}_{x}\\
 & \underline{h}_{y} \ar[u]_{f} \ar[ul]
}
$$ in $\mathsf{\mathbf{Ho}}(\mathsf{\mathbf{Pr}}_{M}(C))$. It is clear that $\mathsf{\mathbf{Ho}}(\underline{R})$ is a (usual) sieve over $x$ in $\mathsf{\mathbf{Ho}}_{M}(C)$.

\begin{prop}\label{inducedtop}
If $\underline{\tau}$ is an enriched homotopy topology on an $M$-enriched category $C$, then the set of all $\mathsf{\mathbf{Ho}}(\underline{R})$ for all enriched $\underline{\tau}$-covering sieves $\underline{R}$ defines a Grothendieck topology $\mathsf{\mathbf{Ho}}(\underline{\tau})$ on $\mathsf{\mathbf{Ho}}_{M}(C)$.
\end{prop}

\textbf{Proof.} Left as an easy exercise. \hfill $\mathbf{\Box}$\\


\begin{Q}\emph{\begin{enumerate}\item Is there a converse construction to $\underline{\tau}\longmapsto \mathsf{\mathbf{Ho}}(\underline{\tau})$ ?
	\item Does an enriched homotopy topology $\underline{\tau}$ on an $M$-enriched category $C$, define an enriched topology on the $\mathbf{Ho}(M)$-enriched category $\mathbf{Ho}(C)$ (with trivial model structures) ?
\end{enumerate}
}
\end{Q}
\end{section}

\begin{section}{Enriched homotopical stacks}

\begin{subsection}{Enriched Bousfield localizations}

\begin{df}
If $C$ and $C'$ are $M$-enriched model categories, an \emph{$M$-enriched Quillen adjunction}  is an $M$-enriched adjunction $(L,R)$, where $L:C\rightarrow C'$ and $R:C'\rightarrow C$, such that the underlying adjunction $(L_{0},R_{0})$, where $L_{0}:C_{0}\rightarrow C'_{0}$ and $R_{0}:C'_{0}\rightarrow C_{0}$ is a Quillen adjunction between the underlying model categories.
\end{df}

\begin{df}
Let $C$ be an $M$-enriched model category and $S$ be a set of maps in $C$ (i.e. morphisms in $C_{0}$).
\begin{itemize}
	\item  An object $x$ in $C$ is $S$\emph{-local over $M$} if it is fibrant and  for any $f:y\rightarrow y'$ in $S$, the induced morphism, $f^{*}:\mathbb{R}\underline{Hom}_{M}(y',x)\rightarrow \mathbb{R}\underline{Hom}_{M}(y,x)$ is an isomorphism in $\mathbf{Ho}(M)$. 
	\item A map $u:x\rightarrow x'$ in $C$ is an $S$\emph{-local equivalence over $M$} if for any $S$-local object over $M$ $y$, the canonical morphism $g^{*}:\mathbb{R}\underline{Hom}_{M}(x',y)\rightarrow \mathbb{R}\underline{Hom}_{M}(x,y)$ is an isomorphism in $\mathbf{Ho}(M)$.
\end{itemize}
\end{df}

\begin{thm}\label{enrloc}
Let $C$ be an $M$-enriched model category and $S$ be a set of maps in $C$ (i.e. morphisms in $C_{0}$). 
\begin{itemize}
	\item $C$ endowed with the classes of $S$-local equivalences (as weak equivalences) and cofibrations, is an $M$-enriched model category, denoted as $\mathbf{L}^{M}(C;S)$.
	\item The identity functor $\mathrm{Id}=:\mathrm{loc}_{S}: C\longrightarrow \mathbf{L}^{M}(C;S)$ has the following properties:
\begin{enumerate}
	\item  $\mathrm{loc}_{S}$ is an $M$-enriched left Quillen functor.
	\item for any $M$-enriched model category $C'$ and any $M$-enriched left Quillen functor $L:C\rightarrow C'$ such that $\mathbb{L}L:\mathbf{Ho}(C)\rightarrow \mathbf{Ho}(C')$ sends $S$ to isomorphisms, there exists a factorization $$\xymatrix{C \ar[r]^-{\mathrm{loc}_{S}} \ar[dr]_{L} & \mathbf{L}^{M}(C;S) \ar[d]\\
	& C'}$$
\end{enumerate}

\end{itemize}
The model category $\mathbf{L}^{M}(C;S)$ will be called the $M$\emph{-enriched left Bousfield localization of $C$ with respect to $S$.}
\end{thm}

The proof of Theorem \ref{enrloc} will be given elsewhere; the idea is to adapt the proof of \cite{hi} to the present enriched context. For the additional set-theoretic hypotheses needed on $C$ and $M$, see Remark \ref{universe} below.\\

\end{subsection}

\begin{subsection}{Enriched homotopical $\check{\textrm{C}}$ech-stacks}

The axioms of an enriched model category imply that the functor $$\underline{Hom}_{\mathsf{\mathbf{Pr}}_{M}(C)}(-,-):\mathsf{\mathbf{Pr}}_{M}(C)^{op}\times \mathsf{\mathbf{Pr}}_{M}(C) \longrightarrow M$$ preserves fibrations and trivial fibrations hence can be derived to the right; therefore, for any $F,F'\in \mathsf{\mathbf{Pr}}_{M}(C)$, we denote by  $$\mathbb{R}\underline{Hom}_{\mathsf{\mathbf{Pr}}_{M}(C)}(F,F'):=\underline{Hom}_{\mathsf{\mathbf{Pr}}_{M}(C)}(QF,RF')$$ the \textit{derived $M$-enriched mapping space} between $F$ and $F'$ which is a well defined object in $\mathsf{\mathbf{Ho}}(M)$.

\begin{df}
Let $(C,\underline{\tau})$ be an $M$-enriched homotopy site. An $M$-valued presheaf $F\in \mathsf{\mathbf{Pr}}_{M}(C)$ on $C$, is called an $M$\emph{-enriched homotopical $\check{\textrm{C}}$ech-stack} on $(C,\underline{\tau})$ if for any $x$ in $C$ and any $(\underline{R}\rightarrow \underline{h}_{x})\in \underline{Cov}_{\underline{\tau}}(x)$, the induced map $$\mathbb{R}\underline{Hom}_{\mathsf{\mathbf{Pr}}_{M}(C)}(\underline{h}_{x},F)\rightarrow \mathbb{R}\underline{Hom}_{\mathsf{\mathbf{Pr}}_{M}(C)}(\underline{R},F)$$
\noindent is an isomorphism in $\mathbf{Ho}(M)$.
\end{df}

\begin{prop}
Let $(C,\underline{\tau})$ be an $M$-enriched homotopy site. There is a model structure on $\mathbf{Pr}_{M}(C)$, denoted as $\check{C}_{M}^{\sim,\underline{\tau}}$, such that $\mathbf{Ho}(\check{C}_{M}^{\sim,\underline{\tau}})$ is equivalent to the full subcategory of $\mathbf{Ho}(\mathbf{Pr}_{M}(C))$ consisting of $M$-enriched homotopical $\check{\textrm{C}}$ech-stacks.
\end{prop}

\textbf{Proof.} Just take $\check{C}_{M}^{\sim,\underline{\tau}}:=\mathbf{L}^{M}(\mathbf{Pr}_{M}(C); S_{\underline{\tau}})$, where $S_{\underline{\tau}}$ is the set of all enriched covering sieves.  \hfill $\mathbf{\Box}$\\

\begin{ex}
\emph{For $M=\mathsf{\mathbf{Set}}$, $\underline{\tau}$ is a usual Grothendieck topology on $C$ and $\mathbf{Ho}(\check{C}_{\mathbf{Set}}^{\sim,\underline{\tau}})\simeq \mathbf{Sh}(C;\underline{\tau})$. For $M=\mathsf{\mathbf{SSet}}$ and $C$ with the trivial $\mathsf{\mathbf{SSet}}$-enrichment, $\check{C}_{\mathbf{SSet}}^{\sim,\underline{\tau}}$ has been considered in \cite{dhi} and \cite{lu}. }
\end{ex}

\end{subsection}

\begin{subsection}{Enriched hyperdescent and enriched homotopical stacks}

\begin{df}\label{coverings}
Let $\underline{\tau}$ be an enriched homotopy topology on an $M$-enriched category $C$. A map $f:F\rightarrow G$ in $\mathsf{\mathbf{Ho}}(\mathsf{\mathbf{Pr}}_{M}(C))$ is a \emph{covering} if for any map $q:\underline{h}_{x}\rightarrow G$ in $\mathsf{\mathbf{Ho}}(\mathsf{\mathbf{Pr}}_{M}(C))$, there exist a $\underline{\tau}$-covering homotopy $M$-sieve $\underline{R}\rightarrow \underline{h}_{x}$, an object $F'$ in $\mathsf{\mathbf{Ho}}(\mathsf{\mathbf{Pr}}_{M}(C))$ and a commutative diagram $$
\xymatrix{ F \ar[rr] & & G\\
F' \ar[u] \ar[r]_{p} & \underline{R} \ar[r] & \underline{h}_{x} \ar[u]_{q}
}$$ where $p$ is a homotopy epimorphism in $\mathsf{\mathbf{Ho}}(\mathsf{\mathbf{Pr}}_{M}(C))$.
\end{df}

We denote by $\mathbf{sPr}_{M}(C)$ the category of simplicial objects in $\mathbf{Pr}_{M}(C)$; this is a simplicial model category (with the Reedy model structure), \textit{tensored} ($F_{*}\in \mathbf{sPr}_{M}(C)$, $K\in \mathbf{SSet}$ $\;\Rightarrow\;$ $\underline{K}\otimes F_{*} \in \mathbf{sPr}_{M}(C)$) and \textit{cotensored} ($F_{*}\in \mathbf{sPr}_{M}(C)$, $K\in \mathbf{SSet}$ $\;\Rightarrow\;$ $F_{*}^{\underline{K}} \in \mathbf{sPr}_{M}(C)$) over $\mathbf{SSets}$. We will simply denote by $F_{*}^{K} \in \mathbf{Pr}_{M}(C)$   
the $0$-th level of the simplicial object $F_{*}^{\underline{K}} \in \mathbf{sPr}_{M}(C)$; in particular, one has a natural isomorphism $F_{*}^{\Delta[n]}\simeq F_{n}$, for any $F_{*} \in \mathbf{sPr}_{M}(C)$, . \\  
Recall that equivalences in $\mathbf{sPr}_{M}(C)$ are levelwise equivalences and the fibrations are the   
morphisms $F_{*} \longrightarrow G_{*}$  
such that, for any $[n] \in \Delta$, the induced morphism  
$$F_{*}^{\Delta[n]} \simeq F_{n} \longrightarrow F_{*}^{\partial \Delta[n]}  
\times_{G_{*}^{\partial \Delta[n]}}G_{*}^{\Delta[n]}$$  
is a fibration in $\mathbf{Pr}_{M}(C)$. \\  
  
For any simplicial set $K \in \mathbf{SSet}$, the functor  
$$\begin{array}{ccc}  
\mathbf{sPr}_{M}(C) & \longrightarrow & \mathbf{Pr}_{M}(C) \\  
 F_{*} & \mapsto & F_{*}^{K}  
\end{array}$$  
is a right Quillen functor whose   
right derived functor will be denoted by  
$$\begin{array}{ccc}  
\mathbf{Ho}(\mathbf{sPr}_{M}(C)) & \longrightarrow & \mathbf{Ho}(\mathbf{Pr}_{M}(C)) \\  
F_{*} & \mapsto & F_{*}^{\mathbb{R}K}.  
\end{array}$$  
  
  
For any $F \in \mathbf{Pr}_{M}(C)$ and $K \in \mathbf{SSet}$, we will simply denote by  
$F^{\mathbb{R}A} \in \mathrm{Ho}(\mathbf{Pr}_{M}(C))$ the object $c(F)_{*}^{\mathbb{R}A}$, where $c(F)_{*}\in \mathbf{sPr}_{M}(C)$ is the constant simplicial object associated to $F$. \\  

\begin{df}\label{d12}  
Let $(C,\underline{\tau})$ be an $M$- enriched model site.  
\begin{enumerate}  
\item A morphism in $\mathbf{sPr}_{M}(C)$  
$$F_{*} \longrightarrow G_{*}$$  
is called a \emph{$M$-enriched} $\underline{\tau}$\emph{-hypercover} if $F_{*}$ is Reedy cofibrant in $\mathbf{sPr}_{M}(C)$ and, for any $n\geq 0$, the induced morphism   
$$F_{*}^{\mathbb{R}\Delta^{n}}\simeq F_{n} \longrightarrow  
F_{*}^{\mathbb{R}\partial \Delta^{n}}\times^{h}_{G_{*}^{\mathbb{R}  
\partial\Delta^{n}}}G_{*}^{\mathbb{R}\Delta^{n}}$$  
is a covering in $\mathbf{Ho}(\mathbf{Pr}_{M}(C))$ according to Definition \ref{coverings}.  
  
\item A morphism in $\mathbf{Ho}(\mathbf{sPr}_{M}(C))$  
$$F_{*} \longrightarrow G_{*}$$  
is called a \emph{$M$-enriched} $\underline{\tau}$\emph{-hypercover} if one of its representatives in $\mathbf{sPr}_{M}(C)$ is.  

\item An $M$-enriched $\underline{\tau}$-hypercover $F_{*} \rightarrow G_{*}$ is called \emph{pseudorepresentable} if, for any $n\geq 0$,   
there exists an isomorphism $F_{n}\simeq \coprod_{i\in I_{n}}\underline{h}_{x_{i}}$ in $\mathbf{Ho}(\mathbf{sPr}_{M}(C))$.

\end{enumerate}  
\end{df}  

\begin{df}\label{hyper}
Let $(C,\underline{\tau})$ be an $M$- enriched model site and let $x$ be an object in $C$.

\begin{itemize}
	\item If $R\underline{h}_{x}$ is a fibrant replacement of $\underline{h}_{x}$ in $\mathbf{Pr}_{M}(C)$, we define $R\underline{h}_{x}^{\Delta[*]}\in \mathbf{sPr}_{M}(C)$ by $$R\underline{h}_{x}^{\Delta[*]}: [n]\longmapsto R\underline{h}_{x}^{\Delta[n]}.$$ 
	\item We denote by $\mathsf{hhc}^{M}(x)$ a set\footnote{This is not a set in general; see Remark \ref{universe} for a correct statement.} of representatives for the isomorphism classes in $$\left\{f: H_{*}\rightarrow R\underline{h}_{x}^{\Delta[*]}| f\; \textsl{is a pseudorepresentable hypercover}\right\}.$$ 

	\item We define the set\footnote{Again, this is not a set in general; see Remark \ref{universe}.} of adjunction morphisms in $\mathbf{Pr}_{M}(C)$ (see \cite[Def. 19.8.1]{hi}) $$\mathsf{HHC}^{M}(C):=\left\{\left|f\right|: \left|H_{*}\right|\rightarrow R\underline{h}_{x}| f\in \mathsf{hhc}^{M}(x),\; x\in C \right\}.$$

\end{itemize}
\end{df}

\begin{df}\label{stacks}
Let $(C,\underline{\tau})$ be an $M$-enriched model site. The \emph{model category of $M$-enriched stacks over} $(C,\underline{\tau})$ is by definition $$C^{\sim, \underline{\tau}}_{M}:= \mathbf{L}^{M}(\mathbf{Pr}_{M}(C), \mathsf{HHC}^{M}(C)).$$
$C^{\sim, \underline{\tau}}_{M}$ is an $M$-enriched model category.
\end{df}

Note that, for any Reedy cofibrant $H_{*}\in \mathbf{sPr}_{M}(C)$, there is a natural isomorphism $\mathrm{hocolim}H_{*}\simeq \left|H_{*}\right|$ in $\mathbf{Ho}(\mathbf{Pr}_{M}(C))$ (\cite[Thm. 19.8.7]{hi}). Therefore, the homotopy category $\mathbf{Ho}(C^{\sim, \underline{\tau}}_{M})$ is equivalent to the full subcategory of $\mathbf{Ho}(\mathbf{Pr}_{M}(C))$ consisting of objects satisfying \textit{$M$-enriched $\underline{\tau}$-hyperdescent}, i.e. such that, for any $x\in C$ and any $(f: H_{*}\rightarrow R\underline{h}_{x}^{\Delta[*]})\in \mathsf{hhc}^{M}(x)$, the natural map  $$F(x)\simeq \mathbb{R}\underline{Hom}_{\mathsf{\mathbf{Pr}}_{M}(C)}(\underline{h}_{x},F)\longrightarrow \mathrm{holim}_{\Delta}\mathbb{R}\underline{Hom}_{\mathsf{\mathbf{Pr}}_{M}(C)}(H_{*},F)$$ is an isomorphism in $\mathbf{Ho}(M)$. In deducing the above hyperdescent condition, we used the enriched Yoneda lemma, and the fact that the proof showing that $\mathrm{Map}_{N}(\mathrm{hocolim}-,-)\simeq \mathrm{holim}\mathrm{Map}_{N}(-,-)$, where $N$ is any model category and $\mathrm{Map}_{N}(-,-)$ is the mapping space in $N$ (see \cite{hi}), can be easily adapted to prove that   $$\mathbb{R}\underline{Hom}_{\mathsf{\mathbf{Pr}}_{M}(C)}(\mathrm{hocolim}-,-)\simeq \mathrm{holim}\mathbb{R}\underline{Hom}_{\mathsf{\mathbf{Pr}}_{M}(C)}(-,-).$$

\begin{rmk}\label{universe}\emph{There are set-theoretic problems to be solved in order the left Bousfield localization process, used in Definition \ref{stacks}, to work properly. We give a brief skecth of how to solve the problem by keeping, as in \cite[Proof of Thm. 3.4.1]{hagI}, only pseudorepresentable hypercovers of a bounded size. I owe this solution to Bertrand To\"en. \\Fix a universe $\mathbb{U}$ (we take \cite[Exp. I \& Appendice]{sga4} as a reference for universes and their use with categories). Suppose $C$ is $\mathbb{U}$-small, that $M$ is a $\mathbb{U}$-cofibrantly generated model category (see \cite[Appendix A.1]{hagI}), and finally that each object $m$ in $M$ is $\alpha_{m}$-small for some $\mathbb{U}$-small cardinal $\alpha_{m}$. Denote by $e(C)$ any fixed $\mathbb{U}$-small cardinal bigger than $\aleph_{0}$ and than $$2^{ \sum_{(x,y)\in \mathrm{Ob}(C)\times \mathrm{Ob}(C)}\alpha_{\underline{Hom}_{C}(x,y)}}.$$ \\ If $F\in\mathbf{Pr}_{M}(C)$ (note that $\mathbf{Pr}_{M}(C)$ is an $M$-enriched $\mathbb{U}$-model category), we say that $F$ is of \textit{$e(C)$-small size} if $$\coprod_{x\in \mathrm{Ob}(C)}F(x)$$ is $e(C)$-small. Now for any $x\in C$, define $\mathsf{hhc}_{e(C)}^{M}(x)$ as for $\mathsf{hhc}^{M}(x)$ but restricting to those pseudorepresentable hypercovers $H_{*}$ over $x$ for which $H_{n}$ is of $e(C)$-small size, for any $n\geq 0$. Note that $\mathsf{hhc}_{e(C)}^{M}(x)$ is a $\mathbb{U}$-small set. Let $\mathsf{HHC}_{e(C)}^{M}(C)$ be constructed as in Definition \ref{hyper}, but starting from $\mathsf{hhc}_{e(C)}^{M}(x)$, $x\in C$. Then $\mathsf{HHC}_{e(C)}^{M}(C)$ is again a $\mathbb{U}$-small set and we can apply the enriched Bousfield localization technique, to correctly define $C^{\sim, \underline{\tau}}_{M}:= \mathbf{L}^{M}(\mathbf{Pr}_{M}(C), \mathsf{HHC}_{e(C)}^{M}(C))$ (see also, for the case of $\mathbf{SSet}$-enrichment, \cite[Thm A.2.4]{hagI}).}
\end{rmk}

Note that we have a natural inclusion $$\mathbf{Ho}(C^{\sim, \underline{\tau}}_{M})\hookrightarrow \mathbf{Ho}(\check{C}^{\sim, \underline{\tau}}_{M})$$ of $M$-enriched homotopy stacks into $M$-enriched homotopy $\check{\textrm{C}}\textrm{ech}$ stacks

\end{subsection}

\end{section}

\begin{section}{Some examples}
\begin{subsection}{Gabriel topologies as enriched topologies}

The definition of enriched homotopy topology gives some interesting objects already in the case where the model structures are everywhere trivial. If $R$ is a ring (not necessarily commutative), P. Gabriel defined in \cite[\S V.2 ]{gab} a notion of idempotent topologizing filter on $R$; we will use the name Gabriel filter here for that and its  equivalent characterization given in \cite[\S VI.4 and \S VI.5]{ste} through axioms T1 to T4; in fact axioms T1-T3 exactly defines a topologizing filter in Gabriel's original definition and then T4 is easily seen to be equivalent to the topologizing filter being idempotent.\\
 
\begin{prop}
Let $R$ be an associative ring with unit and $\mathrm{B}R$ its classifying category (i.e. the category with one object whose endomorphism ring is $R$) viewed as a category enriched over abelian groups. Then there is a bijective correspondence between Gabriel filters (of right ideals) on $R$ and $\mathsf{Ab}$-enriched topologies on $\mathrm{B}R$.
\end{prop}

\textbf{Proof.} It is enough to notice that an $\mathsf{Ab}$-enriched sieve on $\mathrm{B}R$ is exactly a right ideal in $R$ and then to compare the axioms (1)-(4) of an $\mathsf{Ab}$-enriched topology above and at the axioms T1-T4 of a Gabriel filter in \cite[\S VI.4 and \S VI.5]{ste}. \hfill $\mathbf{\Box}$\\

\end{subsection}

\begin{subsection}{Spectra-enriched topologies: the case of $\mathrm{B}A$ where $A$ is a ring spectrum.}
Let $M:=\mathbb{S}-\mathbf{Mod}$ be the symmetric monoidal model category of symmetric spectra with its positive model structure. For any monoid $A$ in $\mathbb{S}-\mathbf{Mod}$ we denote by $A-\mathbf{Alg}$ the model category of left  $A$-algebras. For $A\in \mathbb{S}-\mathbf{Alg}$, we let $C:=\mathrm{B}A$ be the $\mathbb{S}-\mathbf{Mod}$-enriched model category with one
object whose enriched endomorphism is given by the $\mathbb{S}$-module $A$ (and composition given by the product   $\mathbb{S}-\mathbf{Mod}$-morphism $A\otimes_{\mathbb{S}}A\rightarrow A$).\\

The question here is to characterize all the $\mathbb{S}-\mathbf{Mod}$-enriched homotopy topologies on $\mathrm{B}A$.\\
First of all let us observe that there is an equivalence of categories $\mathbf{Pr}_{\mathbb{S}-\mathbf{Mod}}(\mathrm{B}A)\simeq A-\mathbf{Mod}$ (where $A-\mathbf{Mod}$ is the model category of left $A$-modules) which is furthermore a Quillen equivalence. Therefore a $\mathbb{S}-\mathbf{Mod}$-enriched homotopy
sieve $\underline{R}$ on the unique object $*$ of $\mathrm{B}A$  can be identified with the pair consisting of the $A$-module $P:=\underline{R}(*)$ together with a homotopy monomorphism $P\rightarrow A$ in $A-\mathbf{Mod}$; any such pair will be called a (left) \textit{ideal} in $A$. Therefore, the $\mathbb{S}-\mathbf{Mod}$-enriched homotopy
sieves in $\mathrm{B}A$ are exactly the (left) ideals of $A$. Since $A-\mathbf{Mod}$ is a stable model category any homotopy monomorphism in $A-\mathbf{Mod}$ is in fact a weak equivalence; therefore the $\mathbb{S}-\mathbf{Mod}$-enriched homotopy
sieves in $\mathrm{B}A$ are exactly the free $A$-modules of rank $1$ (up to equivalence). As observed in the Introduction, this is the prototypical example showing that enriched homotopy topologies (at least as defined in this note) are not really interesting in the case the enrichment is defined over a \textit{stable} model category (like $\mathbb{S}-\mathbf{Mod}$ here).
\end{subsection}

\begin{subsection}{$M=\mathbf{SSet}$ and $C=T$ an $S$-category. Comparison with $S$-sites and stacks over them.}
We want to compare the notion of enriched topology introduced here, in the case $M=\mathbf{SSet}$, with the notion of $S$-topology introduced in \cite{hagI}.\\
Let $M=\mathbf{SSet}$ and $T$ any $S$-category. Recall that an $S$-topology on $T$ is just a Grothendieck topology on $\mathbf{Ho}(T)$. There is a $0$-truncation map
$$\pi_{0}:\left\{\mathbf{SSet}\textrm{-enriched homotopy topologies on}\; T \right\}\longrightarrow \left\{S\textrm{-topologies on}\; T\right\}$$
that sends a $\mathbf{SSet}$-enriched $\underline{\tau}$-sieve $\underline{R}$ to $\pi_{0}(\underline{R})$ which is a sieve in $\pi_{0}T$. Moreover, we have also a map in the other direction
$$\underline{(-)}:\left\{S\textrm{-topologies on}\; T \right\} \longrightarrow \left\{\mathbf{SSet}\textrm{-enriched homotopy topologies on}\; T \right\}$$
defined as follows. If $i:R\hookrightarrow [-,x]$ is a $\tau$-sieve over $x\in T$, for some $S$-topology $\tau$ on $T$, we define the induced $\underline{\tau}$-sieve $\underline{R}$ through the following homotopy cartesian diagram in $\mathbf{Pr}_{\mathbf{SSet}}(T)=\mathbf{SPr}(T)$

$$\xymatrix{ \underline{R} \ar[r]^{i'} \ar[d] & \underline{h}_{x} \ar[d] \\ 
             cR \ar[r]_{i}                    & c[-,x] }
$$

\noindent where $c$ denotes the constant simplicial set functor $\mathbf{Set} \rightarrow \mathbf{SSet}$. Note that by Proposition \ref{hpbmono}, $i'$ is still a homotopy monomorphism. It is not difficult to check that this actually gives a $\mathbf{SSet}$-enriched homotopy topology on $T$ and that, if $\tau$ is an $S$-topology, then $\pi_{0}(\underline{\tau})= \tau$.\\

\begin{prop}
Let $T$ be an $S$-category. We denote by $\mathbf{SPr}_{\tau}(T)$ the model category of stacks on the $S$-site $(T,\tau)$, defined in \cite[\S 3.4]{hagI}. 
\begin{itemize}
	\item If $\tau$ is an $S$-topology on $T$, then there is a canonical equivalence $\mathbf{Ho}(T_{\mathbf{SSet}}^{\underline{\tau},\sim})\simeq \mathbf{Ho}(\mathbf{SPr}_{\tau}(T)).$
	\item If $\underline{\tau}$ is a $\mathbf{SSet}$-enriched topology on $T$, then there is a canonical equivalence  $\mathbf{Ho}(T_{\mathbf{SSet}}^{\underline{\tau},\sim})\simeq \mathbf{Ho}(\mathbf{SPr}_{\pi_{0}(\underline{\tau})}(T).$
\end{itemize}
\end{prop}

\end{subsection}
\begin{subsection}{DG-categories}\label{dgcats}
DG-categories $C\equiv\mathcal{A}$ fit in our picture with $M=\mathbf{Ch}(k)$ (the unbounded
category of complexes of modules over some $\mathbb{Q}$-algebra $k$ with its projective model structure). Note that $\mathbf{Pr}_{M}(\mathcal{A})\simeq \mathcal{A}\mathbf{-Mod}$, see e.g. \cite[App. III]{dr}. As in the spectra-enriched case, the problem here is that $\mathbf{Ch}(k)$ is a stable model category; therefore our theory does not seem to yield interesting new constructions here.\\ However one could consider the following remedy (that suitably varied should apply also to the spectra-enriched case). Every DG-category $\mathcal{A}$ has an ``underlying'' $S$-category (i.e. a simplicially enriched category) $S(\mathcal{A})$ which is the Dwyer-Kan localization (with respect to weak equivalences) of the category of quasi-representable left $\mathcal{A}^{\mathrm{opp}}$-modules (\cite[\S 3]{dg}). Therefore one could decide to look for ``enriched homotopy topologies'' on 
$\mathcal{A}$ as $\mathbf{SSet}$-enriched homotopy topologies on $S(\mathcal{A})$. I have not investigated this point further, but if it works it could also provide a similar way-out from the difficulties occurring in other situations where a stable model category-enrichment is given.
\end{subsection}

\begin{subsection}{Some questions}
Here are a few questions suggested to the interested reader.
\begin{itemize}
	\item What are the enriched homotopy topologies on a 2-category ?
	\item Does one get something more interesting in the context of Section \ref{dgcats} if one takes the enrichment in non-negatively graded chain complexes of $k$-modules ?
	\item As showed above, the theory presented in this note is not very exciting when the enrichment is prescribed over a \textit{stable} model category, the reason being essentially that homotopy monomorphisms (hence enriched homotopy sieves) are exactly weak equivalences in this case.\\ Does the remedy proposed in Section \ref{dgcats} extend to cover the case of enrichments in other stable model categories (e.g. in spectra) ?
\end{itemize}
\end{subsection}
 
\end{section}

\begin{section}{Another approach (without sieves)}
Consider $\underline{h}_{-}:C \rightarrow \mathbf{Pr}_{M}(C)$ and define the category $\mathrm{Im}^{h}(\underline{h}_{-})$ as the full subcategory of $\mathbf{Pr}_{M}(C)$ consisting of objects that are isomorphic in $\mathbf{Ho}(\mathbf{Pr}_{M}(C))$ to objects of the form $\underline{h}_{x}$, for $x\in C$. If we suppose that the $M$-enrichment is such that $\mathrm{Im}^{h}(\underline{h}_{-})$ is closed under homotopy pullbacks computed in the model category $\mathbf{Pr}_{M}(C)$, then $\mathrm{Im}^{h}(\underline{h}_{-})$ is a \textit{pseudo-model category} (\cite[\S 4]{hagI}), and we denote it by $C^{\wedge}_{\mathrm{repr},\, M}$. Then we know what a model pretopology on $C^{\wedge}_{\mathrm{repr},\, M}$ is (\cite[\S 4]{hagI}) and we might wish to define an $M$-enriched model pretopology on $C$ to be just a model pretopology on the pseudo-model category $C^{\wedge}_{\mathrm{repr},\, M}$. This would give, taking into account the $M$-enrichment of $C^{\wedge}_{\mathrm{repr},\, M}$, an associated theory of ``$M$-enriched homotopy stacks'', along the lines of \cite[\S 4]{hagI}.\\
The first question here is to compare this notions to the ones given in \cite[\S 4]{hagI} (and by comparison, above in \S 4.3), in the case $M=\mathbf{SSet}$.\\
A further question is whether this new notion is more suitable to treat the case where $M$ is a stable model category than the previous one (i.e the one using homotopy sieves).
\end{section}

\end{document}